\documentclass[leqno,10pt]{amsart}

\usepackage{amsfonts}
\usepackage{amsmath}
\usepackage{amsfonts}
\usepackage{amssymb}
\usepackage{amscd}
\newtheorem{theorem}{Theorem}
\newtheorem*{acknowledgement*}{Acknowledgement}

\newtheorem*{example*}{Example}

\newtheorem{lemma}[theorem]{Lemma}

\newtheorem{remark}[theorem]{Remark}

\newcommand{\Ric}{\operatorname{Ric}}

\newcommand{\diam}{\operatorname{diam}}
\newcommand{\conj}{\operatorname{conj}}
\newcommand{\Vol}{\operatorname{Vol}}
\newcommand{\ind}{\operatorname{ind}}
\newcommand{\tu}{\tilde{u}}
\newcommand{\tue}{\tilde{u}_{\epsilon}}

\begin{document}

\title[Myers' type theorems]{Myers' type theorems and some related oscillation results}

\author{Paolo Mastrolia}
\address{Dipartimento di Matematica\\
Universit\`a degli Studi di Milano\\
via Saldini 50\\
I-20133 Milano, ITALY} \email{paolo.mastrolia@unimi.it}

\author{Michele Rimoldi}
\address{Dipartimento di Matematica\\
Universit\`a degli Studi di Milano\\
via Saldini 50\\
I-20133 Milano, ITALY} \email{michele.rimoldi@unimi.it}

\author{Giona Veronelli}
\address{Dipartimento di Matematica\\
Universit\`a degli Studi di Milano\\
via Saldini 50\\
I-20133 Milano, ITALY} \email{giona.veronelli@unimi.it}

\date{\today}


\begin{abstract}
In this paper we study the behavior of solutions of a second order
differential equation. The existence of a zero and its
localization allow us to get some compactness results. In
particular we obtain a Myers' type
theorem even in the presence of an amount of negative curvature. The
technique we use also applies to the study of spectral properties of Schr\"odinger
operators on complete manifolds.
\end{abstract}

\maketitle

\section{Introduction and main results}

In 1941 S. B. Myers, \cite{Myers}, obtained his well known and celebrated compactness theorem stating that a
complete $m$-dimensional Riemannian manifold $M$ is compact provided
its Ricci curvature is bounded from below by a positive constant. By the way, Myers proof also permits to get sharp upper diameter estimates. Since
then, this result has been widely extended and improved in several directions. For example, G.J. Galloway,
\cite{Gal-JDG}, proved compactness and a diameter estimate for $M$ perturbing the constant lower bound for the
Ricci curvature by the derivative in radial direction of some bounded function. This is relevant e.g. in the (elliptic) Ricci solitons theory, \cite{PRRS}.

\begin{theorem}[Myers-Galloway]\label{Galloway-Myers}
Let $M$ be a complete Riemannian manifold. Given two different
points $p,q\in M$, let $\gamma_{p,q}$ be a minimizing geodesic
from $p$ to $q$ parameterized by arc length. Suppose that there
exist constants $c>0$ and $F\geq0$ such that for each pair of
points $p,q$ it holds
\[
\Ric(\dot{\gamma}_{p,q},\dot{\gamma}_{p,q})|_{\gamma_{p,q}(t)}\geq c+\frac{d}{dt}\left(f\circ\gamma_{p,q}\right),
\]
for some $C^1(M)$ function $f$ satisfing $\sup_M|f|\leq F$. Then $M$ is compact and
\begin{equation}\label{estimate_diam}
\operatorname{diam}(M)\leq \frac{1}{c}\left[2F+\sqrt{4F^2+\pi^2(m-1)c}\right].
\end{equation}
\end{theorem}
Myers' proof (and Galloway's generalization) is based on the fact
that, by the second variation formula for the lenght functional,
given a minimizing geodesic $\gamma(s):[0,a]\to M$ between two
points and a smooth function $u$ satisfying $u(0)=u(a)=0$, it
holds
\begin{equation}\label{Myers_var}
0\leq\int_0^a\left[\left(\frac{du}{ds}\right)^2-\frac{\Ric(\dot{\gamma},\dot{\gamma})}{m-1}u^2(s)\right]ds.
\end{equation}
Therefore, if \eqref{Myers_var} is not satisfied for a suitable
choice of $u$ it follows that $\gamma$ is not minimizing and
repeating the argument for each $\gamma$ gives the desired
conclusion. By the way, as first pointed out by G.J. Galloway in
\cite{Gal-PAMS}, the validity of \eqref{Myers_var} and an
integration by parts show that the compactness of $M$ depends on
the behavior, and in particular on the position of the zeros, of
the solution of the differential equation along minimizing geodesics
\begin{equation}\label{main_eq}
Ju(t)=0,
\end{equation}
where the differential operator $J$ is defined as
\[
Ju(t):=-u''(t)-\frac{\Ric (\dot \gamma,\dot\gamma)}{m-1}u(t).
\]
See also \cite{Cal}.

Thus we are reduced to find
sufficient conditions on the Ricci curvature for which the
solutions of the differential equation (\ref{main_eq}) have a
first zero at finite time. At this point, usually one applies
oscillation theory to get geometric assumptions to guarantee that
$M$ is compact; we refer to \cite{Swa} and \cite{Gal-PAMS} for a
more detailed discussion on oscillation theory and compactness.
In particular, using a result by R. Moore, see \cite{Moo}, we get the
following theorem. The case $\lambda=0$ was previously obtained by W. Ambrose in \cite{Amb}. In what follows we denote
\[
\ K_{\gamma}=\frac{\Ric(\dot\gamma,\dot\gamma)}{m-1}.
\]
\begin{theorem}[Ambrose-Moore]\label{Moore}
Suppose there is a point $q\in M$ such that along each geodesic
$\gamma:[0,+\infty)\to M$ parameterized by arc length with
$\gamma(0)=q$ the condition
\begin{equation}\label{eq_ambrose}
\int_0^{\infty}t^{\lambda}K_{\gamma}(t)dt=+\infty
\end{equation}
holds for some $0 \leq \lambda < 1$. Then $M$ is compact.
\end{theorem}

Under the further assumption $\Ric\geq 0$, condition \eqref{eq_ambrose} can be improved. The following result applies a Nehari's oscillation theorem, see \cite{Nea}.

\begin{theorem}[Nehari]\label{Nehari}
Let $\Ric\geq 0$. Suppose there is a point $q\in M$ such that
along each geodesic $\gamma:[0,+\infty)\to M$ parameterized by arc
length with $\gamma(0)=q$ the condition
\[
\int_{t_0}^{\infty}t^{\lambda}K_{\gamma}(t)dt > \frac{(2-\lambda)^2}{4(1-\lambda)}\frac{1}{t_0^{1-\lambda}}
\]
holds for some $t_0>0$ and $0 \leq \lambda < 1$. Then $M$ is compact.
\end{theorem}

As a matter of fact, as we observed above, to conclude that $M$ is compact oscillation theory is not strictly necessary and one could improve Theorem \ref{Moore} and Theorem \ref{Nehari} by focusing his attention upon the more general problem of the existence of a zero for solutions of \eqref{main_eq}. To the best of our knowledge, a few steps have been done in this direction. We point out the paper \cite{Cal} by E. Calabi where the same conclusion of Theorem \ref{Nehari} is reached under assumptions which seem to be neither weaker nor stronger than those of Nehari's result.

\begin{theorem}[Calabi]\label{Calabi}
Let $\Ric\geq 0$. Suppose there is a point $q\in M$ such that
along each geodesic $\gamma:[0,+\infty)\to M$ parameterized by arc
length with $\gamma(0)=q$ it holds
\[
\limsup_{a\to+\infty}\left\{\int_{0}^{a}\sqrt{K_{\gamma}(t)}dt-\frac{1}{2\sqrt{m-1}}\log a\right\}=+\infty.
\]
Then $M$ is compact.
\end{theorem}

Adapting the techniques introduced by Calabi, we are able to extend Theorem
\ref{Moore} and Theorem \ref{Nehari} to the case where the Ricci tensor is
bounded from below by a negative constant. Namely, we obtain the
following result.

\begin{theorem}\label{th_main-B2}
Let $\Ric\geq -(m-1)B^2$, for some constant $B\geq0$. Suppose
there is a point $q\in M$ such that along each geodesic
$\gamma:[0,+\infty)\to M$ parameterized by arc length, with
$\gamma(0)=q$, it holds either
\begin{align}\label{ass_main-B2_0}
    \int_a^btK_{\gamma}(t)dt
    >
    B\left\{b+a\frac{e^{2Ba}+1}{e^{2Ba}-1}\right\}+\frac{1}{4}\log\left(\frac{b}{a}\right).
\end{align}
or
\begin{equation}\label{ass_main-B2}
\int_a^bt^{\lambda}K_{\gamma}(t)dt
>
B\left\{b^{\lambda}+a^{\lambda}\frac{e^{2Ba}+1}{e^{2Ba}-1}\right\}+\frac{\lambda^2}{4(1-\lambda)}\left\{a^{\lambda-1}-b^{\lambda-1}\right\}
\end{equation}
for some $0<a<b$ and $\lambda\neq1$. Then $M$ is compact.
\end{theorem}

\begin{remark}
{\rm In case $B=0$ the expressions in Theorem \ref{th_main-B2} have to be intended in a limit sense. Namely \eqref{ass_main-B2_0} and \eqref{ass_main-B2} have to be replaced respectively by
\begin{align}\tag{\ref{ass_main-B2_0}'} \int_a^btK_{\gamma}(t)dt>1+\frac{1}{4}\log\left(\frac{b}{a}\right)
\end{align}
and
\begin{align}\tag{\ref{ass_main-B2}'}   \int_a^bt^{\lambda}K_{\gamma}(t)dt>\frac{\left(2-\lambda\right)^2}{4\left(1-\lambda\right)a^{1-\lambda}}-\frac{\lambda^2}{4\left(1-\lambda\right)b^{1-\lambda}}
\end{align}
Moreover we note that for $B>0$ and $\lambda=0$ assumption \eqref{ass_main-B2} has the more compact expression
\begin{align}\tag{\ref{ass_main-B2}''}
    (1-e^{-2Ba})\int_a^bK_{\gamma}(t)dt>2B.
\end{align}}
\end{remark}

\begin{remark}
{\rm
Consider a manifold $M$ and its universal covering $\tilde M$.
Since the projection $\pi_M:\tilde M\to M$ is a local isometry we
note that geodesics of $M$ (not necessarily minimizing) lift to
geodesics of $\tilde M$ and Ricci curvature is preserved.
Supposing we are in the assumptions of one of the theorems above,
we have that also $\tilde M$ satisfies the same set of assumptions
and so it is compact. Hence, as observed in \cite{Gal-JDG}, we can
also conclude that the fundamental group $\pi_1(M)$ is finite.}
\end{remark}

Theorem \ref{th_main-B2} will be proved by finding lower and upper bounds for solutions of (\ref{main_eq}). This in turn permits to localize the zeros, if any. The same technique can be used to study solutions $z(t)$ of the more general equation

\begin{equation}\label{main_eq_W}
\begin{cases}&(v(t)z'(t))'+W(t)v(t)z(t)=0\ \textrm{on }(0,+\infty)\\
&z'(t)=O(1)\ \textrm{as }t\searrow0^+,\qquad z(0^+)=z_0>0.
\end{cases}
\end{equation}

In case $v(t)$ and $W(t)$ are nonnegative functions and satisfy very weak re\-gularity and integrability assumptions, equation \eqref{main_eq_W} has been intensively studied by B. Bianchini, L. Mari and M. Rigoli in \cite{BMR}. In particular they dealt with the problem of the existence of a first zero, they studied conditions which imply oscillation and obtained an estimate on the distance of two subsequent zeros. Here we will study the case where $W(t)$ is not necessarily nonnegative, but satisfies the request $Wv^2\geq-B^2$, for some constant $B$. First, we give an integral assumption on $Wv$ which guarantees the existence of a first zero.

\begin{theorem}\label{th_zero}
Let $v(t)$ and $W(t)$ be $L^{\infty}_{loc}([0,+\infty))$ functions such that
\begin{equation}\label{ass_v}
v(t)\geq0,\qquad v(t)^{-1}\in L^{\infty}_{loc}((0,+\infty)),\qquad v^{-1}\notin L^1(0^+),\qquad\lim_{t\to0^+}v(t)=0
\end{equation}
and
\begin{equation}\label{ass_W}
W(t)\geq -\frac{B^2}{v(t)^2}
\end{equation}
for some real constant $B \geq 0$. Let $z(t)\in
Lip_{loc}([0,+\infty))$ be a solution of problem
\eqref{main_eq_W}. If $z(t)\neq 0$ for all $t\in(0,+\infty)$,
then, defining $V(t_1,t_2):=e^{2B\int_{t_1}^{t_2}\frac{ds}{v(s)}}$
for every $t_1,t_2\in[0,+\infty]$, it holds
\[
\int_a^bW(s)v(s)ds\leq \begin{cases}2B &\textrm{if }v^{-1}\notin L^1(+\infty)\\
2B\frac{V(b,+\infty)}{V(b,+\infty)-1}&\textrm{if }v^{-1}\in L^1(+\infty)\end{cases}
\]
for every $0\leq a < b$.

\end{theorem}
Then, iterating the technique of the proof of Theorem \ref{th_zero}, we get an asymptotic condition providing the oscillatory behavior of $z$.
\begin{theorem}\label{th_osc}
Let $v,W$ and $z$ be defined as in Theorem \ref{th_zero}. Then $z$
is oscillatory provided either $v^{-1}\in L^1(+\infty)$ and
\begin{equation}\label{osc_1}
\limsup_{t\to\infty}\left\{\int_R^tW(s)v(s)ds\int_t^{\infty}\frac{ds}{v(s)}\right\}>1
\end{equation}
for some $R>0$, or $v^{-1}\notin L^1(+\infty)$ and
\begin{equation}\label{limsup_not}
\lim_{t\to\infty}\left\{\sup_{t\leq q_1<q_2\leq\infty}\int_{q_1}^{q_2}W(s)v(s)ds\right\}>2B.
\end{equation}
\end{theorem}

\begin{remark}
{\rm
In \cite{Hil}, E. Hille studied the differential equation
\begin{equation}\label{hille}
u''(t)+f(t)u(t)=0
\end{equation}
with $f$ a nonnegative function. In particular he defined the
function $g(x)=x\int_x^{\infty}f(t)dt$ and showed that if
\eqref{hille} is non-oscillatory, then $\liminf g\leq\frac{1}{4}$
and $\limsup g \leq 1$. In case $v^{-1}\in L^{1}(+\infty)$, R.
Moore, \cite{Moo} adapted the first of these conditions to study
equation (\ref{main_eq_W}), showing that $z$ is oscillatory
provided
\[
\liminf_{t\to\infty}\left\{\int_R^t W(s)v(s)ds\int_t^{\infty}\frac{ds}{v(s)}\right\}\geq c,
\]
for some constant $c>\frac{1}{4}$, without any sign assumption on $W$. Up to imposing \eqref{ass_W}, condition (\ref{osc_1}) of Theorem \ref{th_osc} is, in a sense, a ``$\limsup$ counterpart'' of Moore's result.}
\end{remark}

\begin{remark}
{\rm
When $W\geq 0$ (i.e. for $B=0$) and $v^{-1}\in L^1(+\infty)$, in \cite{BMR} the authors defined a critical function
\[
\chi(t):=\left[\left(-\frac{1}{2}\log\int_t^{\infty}\frac{ds}{v(s)}\right)'\right]^2
\]
and prove that $z$ is oscillatory provided
\begin{align}\label{BMR_osc}
\limsup_{t\to+\infty}\int_{T}^t\left(\sqrt{W(s)}-\sqrt{\chi(s)}\right)ds=+\infty
\end{align}
for some constant $T>0$. An easy computation shows that condition \eqref{BMR_osc} is equivalent to
\begin{align}\label{BMR_osc2}
\limsup_{t\to+\infty}\left\{e^{2\int_{T}^t\sqrt{W(s)}ds}\int_t^{\infty}\frac{ds}{v(s)}\right\}=+\infty.
\end{align}
The relation beetween \eqref{BMR_osc2} and \eqref{osc_1} is not so clear. Apparently condition \eqref{BMR_osc2} does not completely contain assumption \eqref{osc_1}.}
\end{remark}

\begin{remark}
{\rm
If $v^{-1}\notin L^1(+\infty)$ we can deduce that $z$ is
oscillatory provided $\int_R^{+\infty}Wv=+\infty$. This was
obtained by W. Leighton without any sign assumption on $W$; see
\cite{Lei}.}
\end{remark}

Now, consider the Schr\"odinger operator $L_w=\Delta+w(x)$, where
$w\in C^0(M)$ and $\Delta$ is the Laplace-Beltrami operator on a
complete non-compact Riemannian manifold $M$. Denote with $B_t$
the geodesic ball centered at some origin $o\in M$, define
$v(t)=\Vol (\partial B_t)$ and let $W(t)$ be the spherical mean of
the potential $w$, that is
\[
\ W(t)=(\Vol (\partial
B_t))^{-1}\int_{\partial B_t}w(x)d\sigma
\]
integrated in the $(m-1)$-dimensional Hausdorff measure $d\sigma$.  By Rayleigh characterization the bottom of the spectrum on the bounded domain $\Omega\subset M$ is defined as
\[
\ \lambda_1^{L_w}(\Omega):=\inf_{\varphi\in Lip_0(\Omega)}\frac{\int_\Omega\left|\nabla \varphi\right|^2-\int_{\Omega}w\varphi^2}{\int_\Omega\varphi^2}
\]
Suppose there exists a solution $z\in Lip(\Omega)$ of the problem \eqref{main_eq_W} with $z(0)=z(t_2)=0$ for some $t_2>0$. As in \cite{BMR}, we consider the function $\varphi_z:M\to \mathbb{R}$ defined as
\[
\ \varphi_z(x):=\begin{cases}z(r(x))&r(x)\leq t_2\\
0&r(x)>t_2,
\end{cases}
\]
where $r(x)$ is the distance function from $o\in M$.
Then, integrating by parts,
\begin{align*}
\lambda_1^{L_w}(B_{t_3})&\leq\frac{\int_{B_{t_3}}\left|\nabla \varphi_z\right|^2-\int_{B_{t_3}}w\varphi_z^2}{\int_{B_{t_3}}\varphi_z^2}=\frac{\int_0^{t_3}v(z^\prime)^2-\int_0^{t_3}Wvz^2}{\int_{0}^{t_3}vz^2}\\&=-\frac{\int_0^{t_3}\left[\left(vz^\prime\right)^{\prime}+Wvz\right]z}{\int_{0}^{t_3}vz^2}=0,
\end{align*}
for every $t_3>t_2$. By the domain monotonicity of eigenvalues we get $\lambda_1^{L_w}(M)<0$.
Moreover, with analogous computations, one can prove that the
oscillation of solutions imply
\begin{align}\label{lambdam-r}
\lambda_1^{L_w}(M\setminus B_R)<0,\qquad\textrm{for all }R\geq 0.
\end{align}
Recall that, given a bounded domain $\Omega\subset M$, the index of $L_w$ is defined as the number of negative eigenvalues of $-L_w$. Hence conditions (\ref{lambdam-r}) together with a result by D. Fisher-Colbrie, \cite{FC}, gives that $L_w$ has infinite index, that is
\[
\ind_{L_w}(M):=\sup_{\Omega\subset M\ \textrm{bounded}}\left\{\ind_{L_w}(\Omega)\right\}=+\infty.
\]
This is the content of the next

\begin{theorem}\label{spectral}
Let $w(x)\in C^0(M)$ be defined on a complete non-compact Riemannian manifold $M$. Let $B\geq0$ be a constant and set $v(r):=\Vol (\partial B_r)$. Suppose that the spherical mean $W(r)$ of $w(x)$ satisfies
\[
W(r)v^2(r)\geq -B^2
\]
for all $r>0$.
\begin{itemize}
\item[i)]Defining $V(t_1,t_2)$ as in Theorem \ref{th_zero}, then $\lambda_1^{L_w}(M)<0$ provided there exist $0<a<b$ such that
\[
\int_{B_{b}\setminus B_{a}}w(x)dx > \begin{cases}2B &\textrm{if }v^{-1}\notin L^1(+\infty)\\
2B\frac{V(b,+\infty)}{V(b,+\infty)-1}&\textrm{if }v^{-1}\in L^1(+\infty).\end{cases}
\]
\item[ii)] $L_w$ is unstable at infinity, i.e. $\lambda_1^{L_w}(M\setminus B_R)<0$ for every $R>0$, provided either $v\in L^1(+\infty)$ and
\begin{equation*}
\limsup_{t\to\infty}\left\{\int_{B_t\setminus B_R}w(x)dx\int_t^{\infty}\frac{ds}{v(s)}\right\}>1
\end{equation*}
for some $R>0$, or $v\notin L^1(+\infty)$ and
\begin{equation*}
\lim_{t\to\infty}\left\{\sup_{t\leq q_1<q_2\leq\infty}\int_{B_{q_2}\setminus B_{q_1}}w(x)dx\right\}>2B.
\end{equation*}
\end{itemize}
In particular, in the assumption \textrm{ii)} $L_w$ has infinite
index.
\end{theorem}

As observed in \cite{BMR}, if $s\left(x\right)$ denotes the scalar curvature of the $m$-dimensional manifold $(M, \left\langle ,\right\rangle)$ and $c_m=4(m-1)/(m-2)$, then, setting $w(x)=-c_m^{-1}s(x)$, the negativity of $\lambda_1^{L_w}$ can be used to prove the existence of positive solutions $u$ of the Yamabe equation
\[
c_m\Delta u + s(x)u - k(x)u^{\frac{m+2}{m-2}}=0.
\]
Here $k(x)$ is the prescribed scalar
curvature of the conformally deformed metric $\widetilde{\left\langle ,\right\rangle}=u^{\frac{4}{m-2}}\left\langle ,\right\rangle$. Hence we obtain the following

\begin{theorem}\label{yamabe}
Suppose that the dimension of $M$ is $m\geq3$ and that the spherical mean $S(r)$ of $s(x)$
satisfies
\[
S(r)\leq c_m\frac{B^2}{\Vol (\partial B_r)}, \qquad r>0,
\]
for some positive constant $B$. Let $k(x)\in C^{\infty}(M)$ be non-positive on $M$ and strictly negative outside a compact set. Set $\mathcal{K}_0 = k^{-1}(0)$ and
\[
\lambda_1^{L_w}(\mathcal{K}_0) = \sup_D\lambda_1^{L_w}(D),
\]
where
\[
L_w =\Delta-\frac{1}{c_m}
s(x),
\]
and $D$ varies among all open sets with smooth boundary containing $\mathcal K_0$. Suppose
\[
\lambda_1^{L_w}(\mathcal{K}_0) > 0.
\]
Defining $V(t_1,t_2)$ as in Theorem \ref{th_zero}, then the background metric can be conformally deformed to a new metric of scalar curvature $k(x)$ provided there exist $0<a<b$ such that
\[
\int_{B_{b}\setminus B_{a}}(-s(x))dx > \begin{cases}2c_mB &\textrm{if }v^{-1}\notin L^1(+\infty)\\
2c_mB\frac{V(b,+\infty)}{V(b,+\infty)-1}&\textrm{if }v^{-1}\in L^1(+\infty).\end{cases}
\]
\end{theorem}

\section{Compactness results}

Since it will be used in the sequel, observe that the existence of a solution of the Cauchy problems involved in our study is guaranteed by minor changes to Proposition A.1 in \cite{BMR}. In fact both problem (\ref{main_eq_W}) with assumptions on the functions as in Theorem \ref{th_zero} and problem \eqref{main_eq} with initial condition $u(0)=0$ admit a locally Lipschitz solution globally defined in $[0,+\infty)$.\\
To start with, we recall the following well known lemma. For a proof avoiding the use of the second variation formula for arc-length see \cite{PRRS}.

\begin{lemma}\label{main_lem}
Let $\left(M,\left\langle ,\right\rangle\right)$ be a complete Riemannian manifold. Fix $o\in M$ and let $r\left(x\right)=dist\left(x,o\right)$. For any point $q \in M$, let $\gamma_{q}:\left[0,r\left(q\right)\right]\rightarrow M$ be a minimizing geodesic from $o$ to $q$ such that $\left|\dot{\gamma}_{q}\right|=1$. If $g\in Lip_{loc}\left(\mathbb{R}\right)$ is such that $g\left(0\right)=g\left(r\left(q\right)\right)=0$, then  for every $q\in M$, it holds
\[
\ 0\leq \int^{r\left(q\right)}_{0}\left(g^{\prime}\right)^{2}ds-\int^{r\left(q\right)}_{0}g^{2}K_{\gamma_q}(t)ds.
\]
\end{lemma}

Myers' theorem shows how a positive lower bound on the Ricci curvature of $M$ sufficies to conclude that $M$ is compact. Nevertheless Lemma \ref{main_lem} can be used to find weaker conditions for compactness. This is
the content of the next theorem, due to G. J. Galloway,\cite{Gal-PAMS}.

\begin{theorem}[Galloway]\label{Gal2}
Let $M$ be an $m$-dimensional complete Riemannian manifold. Suppose
there exists a point $q\in M$ such that for all geodesic
$\gamma:[0,+\infty)\to M$, parameterized by arc length, with
$\gamma(0)=q$, the differential equation
\begin{equation*}
Ju(t)=-u''(t)-K_{\gamma}(t)u(t)=0
\end{equation*}
has a non trivial weak solution $\tilde{u}$ with
$\tilde{u}(t_1)=\tilde{u}(t_2)=0$ for some $0\leq t_1<t_2$
depending on $\gamma$. Then $M$ is compact and
\begin{equation}\label{gal_diam}
\diam M\leq2\max_{\gamma:\gamma(0)=q}t_2.
\end{equation}
\end{theorem}
For the sake of completeness we provide a somewhat direct proof.
\begin{proof}
First, we fix $\gamma$ and show that $\gamma$ stops minimizing beyond $t_2$. Without loss of generality we can suppose $\gamma$ minimizes distances on $[0,t_2]$. Moreover we can assume $t_2$ is the first zero of $\tu$ greater than $t_1$. This is well defined since $\tu(t)>0$ on $[t_1,t_1+\eta]$ for some $\eta$ small enough. Indeed $\tu$ is an eigenfunction of $J$ on $[t_1,t_2]$ corresponding to the eigenvalue $0$. If, by contradiction $\tu(t)=0$ on a sequence $\left\{t_1+\eta_n\right\}_1^{\infty}$ for some $\eta_n\searrow 0$, it would be $\tu\equiv0$ on $[t_1,t_1+\eta]$ by the unique continuation principle of eigenfuctions. Hence, up to change sign, we take $\tilde{u}> 0$ on $(t_1,t_2)$. Denote the bottom of the spectrum of the operator $J$ restricted to the interval $[t_1,t_2]$ by
\[
\lambda_1^{-J}([t_1,t_2])=\inf_{\begin{array}{ll}\scriptstyle{u\in H^2([t_1,t_2])}\\\scriptstyle{u(t_1)=u(t_2)=0}\end{array}}\frac{\int_{t_1}^{t_2}uJu}{\int_{t_1}^{t_2}u^2}.
\]
For considerations above, it is $\lambda_1^{-J}([t_1,t_2])\leq 0$. On the other hand, by Lemma \ref{main_lem} and integrating by parts, we have that
\begin{align}\label{lam}
    \int_{t_1}^{t_2}u(t)Ju(t)dt= -\int_{t_1}^{t_2}u^2(t)K_{\gamma}(t)dt+\int_{t_1}^{t_2}\left.u'\right.^2(t)dt\geq 0
\end{align}
for all $0\leq u\in Lip_{loc}(\mathbb{R})$ such that
$u(t_{1})=u(t_{2})=0$. In particular, replacing $\tilde u$ to $u$
in \eqref{lam} gives that $\lambda_1^{-J}([t_1,t_2])\geq0$. Thus
$\lambda_1^{-J}([t_1,t_2])=0$. Now, fix $\varepsilon>0$ and define a
new function $\tue$ on $[t_1,t_2+\varepsilon]$ as
\[
\tue(t):=\left\{\begin{array}{ll}\tu(t)  & t\in[t_1,t_2] \\ 0 & t\in [t_2,t_2+\varepsilon]. \end{array}\right.
\]
We have that $\tue\in H^2([t_1,t_2+\varepsilon])$ since it is
$H^2$ on both $[t_1,t_2]$ and $[t_2,t_2+\varepsilon]$ and it is
$Lip_{loc}([t_1,t_2+\varepsilon])$. This gives
\begin{equation}\label{lambdaeps}
\lambda_1^{-J}([t_1,t_2+\varepsilon])=\inf_{\begin{array}{ll}\scriptstyle{u\in H^2([t_1,t_2+\varepsilon])}\\\scriptstyle{u(t_1)=u(t_2+\varepsilon)=0}\end{array}}\frac{\int_{t_1}^{t_2+\varepsilon}uJu}{\int_{t_1}^{t_2+\varepsilon}u^2}\leq\frac{\int_{t_1}^{t_2+\varepsilon}\tilde u_{\varepsilon}J\tilde u_{\varepsilon}}{\int_{t_1}^{t_2+\varepsilon}\tilde u_{\varepsilon}^2}=0.
\end{equation}
We show that the inequality is strict. By contradiction, let $\lambda_1^{-J}([t_1,t_2+\varepsilon])=0$. Since $\tue$ realizes the minimum in \eqref{lambdaeps}, it would be an eigenfunction. Then it would be $\tue\equiv 0$ by unique continuation. Contradiction.\\
Thus there exists an eigenfunction $v$ on $[t_1,t_2+\varepsilon]$ such that $v(t_1)=v(t_2+\varepsilon)=0$, $v\geq0$ and $Jv=\lambda_1^{-J}([t_1,t_2+\varepsilon])v$ is non-positive and not identically $0$. Applying Lemma \ref{main_lem}, we obtain that $\gamma$ can not minimize distances on $[t_1,t_2+\varepsilon]$, hence it stops minimizing at $t_2$ as claimed.\\
Now, fix a point $q\in M$ and let $\Gamma$ be the set of geodesics
$\gamma_q$ parameterized by arc length such that $\gamma_q(0)=q$,
define
\[
\conj(q,\gamma_q):=\inf_{\gamma_q\in\Gamma}\left\{t:\gamma_q\textrm{ does not minimize on }[0,t]\right\}.
\]
Set $\conj(q)=\cup_{\gamma:\gamma(0)=q}\conj(q,\gamma)$. Since $M$ is complete, $M$ is compact provided $\conj(q)$ is bounded(see \cite{Amb}). This is trivial since the function $\conj(q,\gamma)$ is continuous with respect to the outgoing geodesic vector $\dot{\gamma}(0)\in\mathbb{S}^m$ by a result of Morse (Lemma 13.1 in \cite{Mor}).\\
Finally let $p_1,p_2\in M$ and consider the geodesics $\gamma_1$ and $\gamma_2$ joining respectively $p_1$ and $p_2$ to $q$. Both $\gamma_1$ and $\gamma_2$ are shorter than $\max_{\gamma:\gamma(0)=q}t_2$. Hence (\ref{gal_diam}) is proved because of the arbitrarity of $p_1$ and $p_2$.
\end{proof}

In the following proofs we will use a comparison result for Riccati equations, which is a generalization of Corollary 2.2 in \cite{PRS}.

\begin{lemma}[Riccati Comparison]\label{Riccati}
Let $G$ and $0<v$ be $C^0([0,+\infty))$ functions and let $q_i\in AC((\bar t, T_i))$, $i=1,2$, be solutions of the Riccati differential inequalities
\begin{align}\label{Riccati_eq}
    q_1'(t)-\frac{q_1^2(t)}{v(t)}-G(t)\geq 0,\qquad q_2'(t)-\frac{q_2^2(t)}{v(t)}-G(t)\leq 0,
\end{align}
a.e. in $(\bar t,T_i)$ satisfying $q_1(\bar t)=q_2(\bar t)$ for some $\bar{t}>0$. Then $T_1\leq T_2$ and $q_1(t)\geq q_2(t)$ in $[\bar t,T_1)$.\\
Conversely, if $q_i\in AC((T_i, \bar t))$, $i=1,2$, are solutions
of (\ref{Riccati_eq}) a.e. in $(T_i,\bar t)$ satisfying $q_1(\bar
t)=q_2(\bar t)$, then $T_1\geq T_2$ and $q_1(t)\leq q_2(t)$ in
$(T_1,\bar t]$.
\end{lemma}

This lemma is proven with minor changes to the proof of Corollary 2.2 in \cite{PRS} and we refer to this for more details.

\begin{proof}
Let $q_i\in AC((\bar t, T_i))$, $i=1,2$, be solutions of (\ref{Riccati_eq}) a.e. in $(\bar t,T_i)$, with $q_1(\bar t)=q_2(\bar t)$. Setting $y_{i}=-q_{i}$ we obtain that
\begin{align}
    y_1'(t)+\frac{y_1^2(t)}{v(t)}+G(t)\leq 0,\qquad y_2'(t)+\frac{y_2^2(t)}{v(t)}+G(t)\geq 0.
\end{align}
Let $\phi_{i}\in C^{1}\left(\left[\bar{t},T_{i}\right)\right)$ be the positive function on $\left[\bar{t},T_{i}\right)$ defined by
\begin{equation}
\phi_{i}=exp\left\{\int_{\bar{t}}^{t}\left(\frac{y_{i}\left(s\right)}{v\left(s\right)}\right)ds\right\}.
\end{equation}
Then $\phi_{i}\left(\bar{t}\right)=1$, $\phi_{i}>0$ on $\left(\bar{t},T_{i}\right)$, $\phi_{i}^{\prime}\in AC\left(\bar{t},T_{i}\right)$ and a straightforward computation shows that
\begin{align*}
&\phi_{i}^{\prime}\left(t\right)=\frac{y_{i}}{v}\phi_{i}\left(t\right),\\
&\phi_{1}^{\prime}\left(\bar{t}\right)=\frac{y_{1}\left(\bar{t}\right)}{v\left(\bar{t}\right)}\phi_{1}\left(\bar{t}\right)=\frac{y_{2}\left(\bar{t}\right)}{v\left(\bar{t}\right)}\phi_{2}\left(\bar{t}\right)=\phi_{2}^{\prime}\left(\bar{t}\right)
\end{align*}
and
\begin{align}\label{Sturm_eq}
    \left(v\phi_{1}^{\prime}\right)^{\prime}+G\phi_{1}\leq0 \textrm{ a.e. in} \left(\bar{t},T_{1}\right),\qquad \left(v\phi_{2}^{\prime}\right)^{\prime}+G\phi_{2}\leq0 \textrm{ a.e. in} \left(\bar{t},T_{2}\right).
\end{align}
Adapting the Sturm comparison result of Theorem 2.1 in \cite{PRS} to the differential inequalities (\ref{Sturm_eq}) we have that if $\phi_{i}\in C^{1}\left(\left[\bar{t},T_{i}\right)\right)$ are solutions of (\ref{Sturm_eq}) with the properties obtained above then
\begin{align*}
\frac{\phi_{1}^{\prime}}{\phi_{1}}\leq \frac{\phi_{2}^{\prime}}{\phi_{2}}, \qquad T_{1}\leq T_{2}\quad\textrm{and }\qquad \phi_{1}\leq\phi_{2} \textrm{ on } \left[\bar{t}, T_{1}\right).
\end{align*}
This shows that $-q_{1}=y_{1}=\frac{\phi_{1}^{\prime}}{\phi_{1}}v\leq \frac{\phi_{2}^{\prime}}{\phi_{2}}v=y_{2}=-q_{2}$ on $\left(\bar{t},T_{1}\right)$, as required.\\
The second part of the lemma can be proven similarly making a change of variable from $t$ to $-t$.
\end{proof}

We are now in the position to prove Theorem \ref{th_main-B2}.

\begin{proof}[Proof.(of Theorem \ref{th_main-B2})]
First consider the case $B>0$. Suppose $M$ is non compact. By Theorem \ref{Gal2} for each $q\in M$ there exists
a geodesic $\gamma$ parameterized by arc length with $\gamma(0)=q$ such that each non trivial $Lip_{loc}$
solution $u$ of the problem
\[
\begin{cases}
u''(t)+K_{\gamma}(t)u(t)=0\\
u(0)=0,
\end{cases}
\]
which exists by the considerations at the beginning of this section, should satisfy $u(t)\neq 0$ for all $t>0$.
Hence the function $h(t):=-\frac{u'(t)}{u(t)}$ is well defined and continuous in $(0,+\infty)$. Moreover, since
$u''=-K_{\gamma}u\in L^{\infty}_{loc}([0,+\infty))$ implies $u'$ is locally Lipschitz, we have that $h$
satisfies the differential equation
\begin{equation}\label{diff_h}
h'(t)=h^2(t)+K_{\gamma}(t).
\end{equation}
We want to prove that
\begin{equation}\label{estimate}
-\frac{e^{2Bt}+1}{e^{2Bt}-1}\leq \frac{h(t)}{B}\leq 1,
\end{equation}
for all $t>0$. To this purpose consider the functions
\[
\tilde{h}_C(t)=B\frac{C+e^{2Bt}}{C-e^{2Bt}},\qquad C\geq1,
\]
which are solutions of the equation
\[
\tilde{h}'(t)=\tilde{h}^2(t)-B^2
\]
and note that for all $t>0$ the lower bound on Ricci yields $h'(t)\geq\tilde{h}_C'(t)$ each time $h(t)=\tilde{h}_C(t)$. Moreover $h'(t),\tilde h_C'(t)\geq 0$ where $|h(t)|\geq B$ and
\begin{align*}
    &\tilde{h}_{C}(t)\to+\infty,\qquad\textrm{as }t\to \left(\log C/(2B)\right)^{-},\ C>1,\\
    &\tilde{h}_{C}(t)\to-\infty,\qquad\textrm{as }t\to \left(\log C/(2B)\right)^{+},\ C\geq 1.
\end{align*}
By contradiction, suppose there is a value $t_1$ for which $h(t_1)=H_1>B$. Then we have that
\[
\tilde{h}_{C_1}(t_1)=H_1=h(t_1),\qquad\textrm{for }C_1=\frac{H_1+B}{H_1-B}e^{2Bt_1}>1.
\]
Applying the first part of Lemma \ref{Riccati} with $q_1=h$, $q_2=\tilde{h}_{C_1}$, $G\equiv-B^2$, $v\equiv 1$ and $\bar t=t_1$, we can conclude that $h(t)\to +\infty$ as $t\to t_0$ for some $0<t_0<\frac{\log C_1}{2B}$. Thus $h$ is not globally defined. Contradiction. Similarly, suppose there is a value $t_2$ for which
\[
h(t_2)=H_2<-B\frac{e^{2Bt_2}+1}{e^{2Bt_2}-1}.
\]
Then we have that
\[
\tilde{h}_{C_2}(t_2)=H_2=h(t_2),\qquad\textrm{for }C_2=\frac{H_2+B}{H_2-B}e^{2Bt_2}>1.
\]
As above, we achieve a contradiction by applying the second part of Lemma \ref{Riccati} with $q_1=h$, $q_2=\tilde{h}_{C_2}$ and $\bar t=t_2$.\\
Now we want to use (\ref{diff_h}) and (\ref{estimate}) to contradict \eqref{ass_main-B2}. Then, for $\lambda\neq1$,
\begin{align}\label{est_lambda}
\int_a^bt^{\lambda}K_{\gamma}(t)dt
&=\int_a^b(t^{\lambda}h'(t)-t^{\lambda}h^2(t))dt\\
\nonumber&=\int_a^b\left[(t^{\lambda}h(t))'-t^{\lambda}\left(h(t)+\frac{\lambda}{2t}\right)^2+\frac{\lambda^2}{4}t^{\lambda-2}\right]dt\\
\nonumber&\leq b^{\lambda}h(b)-a^{\lambda}h(a)+\frac{\lambda^2}{4(\lambda-1)}\left[b^{\lambda-1}-a^{\lambda-1}\right]\\
\nonumber&\leq
B\left\{b^{\lambda}+a^{\lambda}\frac{e^{2Ba}+1}{e^{2Ba}-1}\right\}+\frac{\lambda^2}{4(1-\lambda)}\left\{a^{\lambda-1}-b^{\lambda-1}\right\}
\end{align}
for all $b>a>0$. The case $\lambda=1$ can be treated similarly. Finally observe that the computations above work even if we intend all the expressions in a limit sense as $B\to0$. This concludes the proof.
\end{proof}

\begin{remark}
{\rm
Reasoning as in the proof of Theorem \ref{th_main-B2}, we can even find diameter estimates as follows. Suppose $\diam M>D$. Hence by Theorem \ref{Gal2} there exists a geodesic ray $\bar{\gamma}$, with $\bar{\gamma}(0)=q$, such that $\bar{\gamma}$ is minimizing at least on $(0,D/2)$. With notations as above, we have that $h$ has to be defined and continuous at least on $(0,D/2)$. In analogy with (\ref{estimate}), this fact and Riccati comparison force $h$ to satisfy
\begin{align}\label{est_diam}
    -B\frac{e^{2Bt}+1}{e^{2Bt}-1}\leq h(t)\leq B\frac{e^{2B\left(\frac{D}{2}-t\right)}+1}{e^{2B\left(\frac{D}{2}-t\right)}-1}.
\end{align}
This estimate, together with the fact that $K_{\gamma}=h'-h^2$, leads to obtain integral conditions on $K_{\gamma}$, in the spirit of (\ref{est_lambda}). For instance one can prove that $\diam M\leq D$ provided that
\[
2\int_0^{D/4}t^2K_{\gamma}(t)dt>D.
\]}
\end{remark}

\section{Oscillatory behavior and spectral applications}
In this final section we give the proofs of the results concerning
the behavior of solutions of problem \eqref{main_eq_W} and their
geometrical applications; for further details on the proof of these latters see \cite{BMR}.

\begin{proof}[Proof.(of Theorem \ref{th_zero})]
By assumption, $z(t)\in Lip_{loc}([0,+\infty))$ is a solution of problem \eqref{main_eq_W} such that $z(t)\neq
0$ for all $t\in(0,+\infty)$. Defining the function $y(t):=-\frac{v(t)z'(t)}{z(t)}$, we have that $y$ is well
defined in $(0,+\infty)$, is locally Lipschitz by considerations as in the proof of Theorem \ref{th_main-B2} and
it satisfy the differential equation
\begin{equation}\label{diff_y}
\begin{cases}
y'(t)=\frac{y^2(t)}{v(t)}+W(t)v(t)
\\y(0)=0.
\end{cases}
\end{equation}
First of all assume that $v^{-1}\notin L^1(+\infty)$. Proceeding as in the proof of Theorem \ref{th_main-B2}, we want to prove that
\begin{equation}\label{estimate_W}
-1\leq \frac{y(t)}{B}\leq 1,
\end{equation}
for all $t>0$. To this purpose consider the one-parameter family of functions
\begin{equation}\label{y_tilde}
\tilde{y}_C(t)=B\frac{C+V(1,t)}{C-V(1,t)},\qquad C>0,
\end{equation}
which are solutions of the equation
\[
\tilde{y}'(t)=\frac{\tilde{y}^2(t)-B^2}{v(t)}
\]
and note that for all $t>0$ the lower bound on $W(t)$ yields $y'(t)\geq\tilde{y}_C'(t)$ each time $y(t)=\tilde{y}_C(t)$. Moreover $y'(t),\tilde y_C'(t)\geq 0$ where $|y(t)|\geq B$ and
\begin{equation}\label{limit}
i)\lim_{t\to0^+}\tilde{y}_C(t)=B^+;\quad ii)\lim_{t\to+\infty}\tilde{y}_C(t)=-B^-;\quad iii)\lim_{t\to t_C^{\pm}}\tilde{y}_C(t)=\mp\infty,
\end{equation}
where $t_C$ is such that $\int_1^{t_C}v^{-1}(s)ds=\log C/(2B)$. By contradiction, suppose there are values $t_i$, $i=1,2$, for which $y(t_i)=Y_i$ with a) $Y_1>B$ or b) $Y_2<-B$. Then we have that
\[
\tilde{y}_{C_i}(t_i)=Y_i=y(t_i),\qquad\textrm{for }C_i=\frac{Y_i+B}{Y_i-B}V(1,t_i).
\]
Choose $q_1=y$, $q_2=\tilde{h}_{C_i}$, $G=-B^2/v$ and $\bar
t=t_i$. Applying the first part of Lemma \ref{Riccati} for $i=1$
and the second part for $i=2$, we can conclude that a) yields
$y(t)\to +\infty$ as $t\to t'^-$ for some $t_1<t'<t_{C_1}$ while
b) leads to conclude that $y(t)\to -\infty$ as $t\to t''^+$ for
some $0<t_{C_2}<t''<t_2$. Thus $y$ is not globally defined. This
contradiction implies the validity of \eqref{estimate_W}, which
gives
\begin{align}\label{Wv_zero}
\int_a^bW(s)v(s)ds\leq\int_a^by'(s)ds=y(b)-y(a)\leq 2B.
\end{align}
Now, let $v^{-1}\in L^1(+\infty)$. In this case the limit (\ref{limit}.iii) holds only for $C<V(1,+\infty)$, since otherwise $\tilde y_C$ is well defined all over $(0,+\infty)$. Note that also (\ref{limit}.ii) is satisfied with a different limit, but this has no importance to our purpose. Hence the estimate \eqref{estimate_W} gets modified in
\begin{align}\label{estimate_W_in}
-1\leq \frac{y(t)}{B}\leq \frac{V(t,+\infty)+1}{V(t,+\infty)-1},
\end{align}
which in turn implies
\begin{align*}
\int_a^bW(s)v(s)ds\leq\int_a^by'(s)ds=y(b)-y(a)\leq \frac{2BV(b,+\infty)}{V(b,+\infty)-1}.
\end{align*}
\end{proof}

\begin{proof}[Proof.(of Theorem \ref{th_osc})]
First, we assume $v^{-1}\notin L^1(+\infty)$ and consider the
functions $\tilde y_C$ defined as in \eqref{y_tilde}. By
contradiction, suppose $z$ is not oscillatory. Hence there exists
$T>0$ such that $z$ has no zeros in $(T,+\infty)$, which in turn
implies that the function $y(t)=-\frac{v(t)z'(t)}{z(t)}$ is
globally defined in this interval. As shown in the proof of
Theorem \ref{th_zero}, this forces $y(t)\leq B$ for all $t>T$. In
fact we can prove
\begin{equation}\label{y-}
-B\frac{V(T,t)+1}{V(T,t)-1}\leq y(t)\leq B,\qquad t>T.
\end{equation}
Indeed the RHS of \eqref{y-} is exactly the function $\tilde
y_{\bar C}$ for $\bar C=V(1,T)$. By (\ref{limit}), we get that,
for $C>\bar C$, $\tilde y_C$ is a monotone non decreasing function
with a vertical asymptote in some $t_C>T$. If there exists a point
$t_1>T$ such that \eqref{y-} is not verified in $t_1$, we
contradict the global definition of $y$ in $(T,+\infty)$ by
applying Lemma \ref{Riccati} as in the previous proofs.\\
Finally, as in \eqref{Wv_zero} we get
\begin{align}\label{est_osc}
    \int_{a}^bW(s)v(s)ds\leq\frac{2BV(T,a)}{V(T,a)-1},
\end{align}
for all $b>a>T$. Hence the existence of such a $T$ contradicts \eqref{limsup_not}, since RHS of \eqref{est_osc} tends to $2B$ as $a\to\infty$.\\
Now, let $v^{-1}\in L^{1}(+\infty)$. As above, suppose $z$ has no zeros in $(T,+\infty)$ for some $T>R$. As in \eqref{estimate_W_in} we get
\[
y(t)\leq B\frac{V(t,+\infty)+1}{V(t,+\infty)-1},
\]
since otherwise $y(t)$ is forced to have a vertical asymptote at some finite $t_0>t$. Moreover, reasoning exactly as in the case $v^{-1}\notin L^1(+\infty)$, we get the lower estimate
\[
y(t)\geq -B\frac{V(T,t)+1}{V(T,t)-1},
\]
for $t>T$. This estimates in turn give
\begin{align}\label{Wv_osc}
    \int_{a}^bW(s)v(s)ds\leq B \left\{ \frac{V(b,+\infty)+1}{V(b,+\infty)-1}+\frac{V(T,a)+1}{V(T,a)-1}\right\},
\end{align}
for all $b>a>T$. By assumption \eqref{osc_1} there exists a $\delta>0$ and a sequence $\left\{t_n\right\}_{n=1}^{\infty}\nearrow +\infty$ such that
\begin{align}\label{3delta}
\int_R^{t_n}W(s)v(s)ds\int_{t_n}^{\infty}\frac{ds}{v(s)}>1+3\delta
\end{align}
for all $n\geq 1$. Since $v^{-1}\in L^1(+\infty)$, there exists $N_1\in\mathbb{N}$ such that
\[
\int_R^{T+1}W(s)v(s)ds\int_{t_n}^{\infty}\frac{ds}{v(s)}<\delta
\]
for all $n>N_1$. This latter combined with \eqref{3delta} gives
\begin{align}\label{2delta}
\int_{T+1}^{t_n}W(s)v(s)ds\int_{t_n}^{\infty}\frac{ds}{v(s)}>1+2\delta
\end{align}
for all $n>N_1$. We note that
\[
\varepsilon\sim e^{\varepsilon}-1 \sim 2\frac{e^{\varepsilon}-1}{e^{\varepsilon}+1},\qquad\textrm{as }\varepsilon\to 0^+.
\]
Since $\int_{t_n}^{\infty}v^{-1}\to 0$ as $n\nearrow\infty$, there exists $N_2\in\mathbb{N}$ such that
\begin{align}\label{asym}
2\frac{V(t_n,+\infty)-1}{V(t_n,+\infty)+1}>\frac{\delta+1}{2\delta+1}2B\int_{t_n}^{\infty}\frac{ds}{v(s)}
\end{align}
for all $n>N_2$. Then (\ref{2delta}) and (\ref{asym}) imply
\begin{align}\label{delta}
\frac{V(t_n,+\infty)-1}{V(t_n,+\infty)+1}\int_{T+1}^{t_n}\frac{W(s)v(s)}{B}ds
&\geq \frac{\delta+1}{2\delta+1}\int_{t_n}^{+\infty}\frac{ds}{v(s)}\int_{T+1}^{t_n}W(s)v(s)ds \\
&> 1+\delta\nonumber
\end{align}
for all $n>\max\left\{N_1,N_2\right\}$. Moreover, since $v^{-1}\in L^1(+\infty)$, (\ref{3delta}) gives
\[
\int_{T+1}^{t_n}W(s)v(s)ds=\left(\int_{R}^{t_n}W(s)v(s)ds-\int_{R}^{T+1}W(s)v(s)ds\right)\nearrow+\infty
\]
as $n\nearrow\infty$, which in turn implies there exists $N_3\in\mathbb{N}$ such that
\begin{equation}\label{N3}
\frac{V(R,T+1)-1}{V(R,T+1)+1}\int_{T+1}^{t_n}W(s)v(s)ds>\frac{(1+\delta)(2+\delta)}{\delta}
\end{equation}
for all $n>N_3$. Choose $a=T+1$ and $b=t_n$. Combining \eqref{Wv_osc}, \eqref{delta} and \eqref{N3} we get
\begin{align}
1&\geq\int_{T+1}^{t_n}\textstyle{\frac{W(s)v(s)}{B}ds}\left\{\textstyle{\frac{V(t_n,+\infty)+1}{V(t_n,+\infty)-1}+\frac{V(R,T+1)+1}{V(R,T+1)-1}}\right\}^{-1}\\
&= \left\{\textstyle{\frac{V(t_n,+\infty)+1}{V(t_n,+\infty)-1}\left(\int_{T+1}^{t_n}\frac{W(s)v(s)}{B}ds\right)^{-1}+\frac{V(R,T+1)+1}{V(R,T+1)-1}\left(\int_{T+1}^{t_n}\frac{W(s)v(s)}{B}ds\right)^{-1}}\right\}^{-1}\nonumber\\
&> \left\{\frac{1}{1+\delta}+\frac{\delta}{\left(1+\delta\right)\left(2+\delta\right)}\right\}^{-1}\nonumber\\
&= 1+\frac{\delta}{2}>1\nonumber
\end{align}
for all $n>\max\left\{N_1,N_2,N_3\right\}$. Contradiction.
\end{proof}

\begin{proof}[Proof.(of Theorem \ref{spectral})]
Proposition 1.2 and the considerations at the beginning of Section 2 yield assumptions
\eqref{ass_v} are satisfied and there exists a locally Lipschitz
solution of \eqref{main_eq}. Then Theorem \ref{spectral} is
implied by Theorems \ref{th_zero} and \ref{th_osc} as in the proof
of Theorem 1.4 in \cite{BMR}.
\end{proof}
\bigskip

\begin{acknowledgement*}
The authors are deeply grateful to Stefano Pigola for his guidance and constant encouragement during the
preparation of the manuscript.
\end{acknowledgement*}

\end{document}